\documentclass{patmorin}
\usepackage{pat}
\usepackage[T1]{fontenc}
\usepackage[utf8]{inputenc}
\usepackage{mathtools}
\usepackage{thmtools}
\usepackage{thm-restate}

\usepackage[inline]{enumitem}
\usepackage{soul}
\usepackage[normalem]{ulem}

\crefname{p}{}{}
\creflabelformat{p}{#2(#1)#3}

\usepackage{todonotes}

\usepackage[longnamesfirst,numbers,sort&compress]{natbib}

\newcommand{\F}{\mathcal F}

\title{\MakeUppercase{The Volume Helly Theorem in the plane, colorful version}%
  \thanks{This research was partly funded by NSERC.}}

\author{
Imre B{\'a}r{\'a}ny%
\thanks{Alfr{\'e}d R{\'e}nyi Institute of Mathematics, HUN-REN,
13 Re{\'a}ltanoda Street, Budapest 1053, Hungary, and 
Department of Mathematics, University College London, Gower Street, London, WC1E 6BT, UK.}
\and\qquad 
Bobby Miraftab%
\thanks{School of Computer Science, Carleton University, Ottawa, Canada.}
\and\quad
Leonidas Theocharous%
\thanks{School of Electrical Engineering and Computer Science, University of Ottawa, Ottawa, Canada.}
}

\date{}

\begin{document}
\maketitle

\begin{abstract}
We prove a colorful volume Helly theorem for convex sets in $\R^2$: There is a constant $V>0$ such that if $\F_1,\F_2,\F_3,\F_4$ are finite families of convex sets in $\R^2$ and if $|\bigcap_1^4F_i|\ge V$ for every transversal $F_i\in \F_i,\; (i=1,\ldots,4)$, then $|\bigcap \F_i|\ge 1$ for some $i$. Here $|A|$ is the Lebesgue measure of $A\subset \R^d$. The main ingredient is the following theorem. Let $Q_1,\ldots,Q_4\subset\mathbb R^2$ be convex quadrilaterals of area at most $1$, where of course each $Q_i$ 
is the intersection of 4 halfplanes. Then for every $Q_i$ there is one of these halfplanes $H_i$, say, such that $|\bigcap_1^4 H_i| \le 4096$. 
\end{abstract}

\section{Introduction}

Helly's theorem is a fundamental local-to-global principle for convexity.
If $\mathcal F$ is a finite family of convex sets in $\mathbb R^d$ and every $d+1$ members of $\mathcal F$ have a common point, then all members of $\mathcal F$ have a common point~\cite{Helly1923}.  Many Helly-type theorems replace the conclusion ``non-empty intersection'' by a quantitative condition on the intersection. One important example is the volume Helly theorem of B\'ar\'any, Katchalski, and Pach~\cite{BaranyKatchalskiPach1984}. 

\begin{thm}\label{th:volHelly} For every $d\ge 2$ there is a number $v(d)>0$ such that the following holds. Let $\F$ be a finite family of convex sets in $\R^d$ with $|\F|\ge 2d$. If $|\bigcap \F^*|> v(d)$ for every subfamily $\F^* \subset \F$ of size $2d$, then $|\bigcap \F| > 1$.
\end{thm}

Simple examples show that  $2d$, the size of the subfamily $\F^*$ in the condition, cannot be reduced. Thus, for volume, the relevant Helly number is $2d$, rather than $d+1$ of the classical nonemptiness theorem.

Another direction is the colorful version of Helly's theorem discovered by Lov\'asz and first published by B\'ar\'any~\cite{Barany1982}. Before stating it we fix some notation and terminology. 

Assume $\F_1,\ldots,\F_k$ are finite (and non-empty) families of convex sets in $\mathbb R^d$, $k\ge 2$. A {\sl transversal} of this system is $F_i\in \F_i$ for every $i\in [k]$. The $\F_i$ are often considered as colors in which case a transversal can be called a {\sl colored choice}. Here comes the colorful Helly theorem.

\begin{thm}\label{th:colHelly} If $\F_1,\ldots, \F_{d+1}$ are finite families of convex sets in $\R^d$ and for every transversal $F_i\in \F_i,\; i\in [d+1]$ of this system $\bigcap F_i\ne \emptyset$, then $\bigcap F_i \ne \emptyset$ for some $i\in [d+1]$. Equivalently, if no color class has a point in common, then there is a rainbow choice with empty intersection.
\end{thm}

It is natural to ask for a theorem that is both colorful and quantitative. Results of this kind are known with more than $2d$, namely $3d$, color classes~\cite{DamasdiFoldvariNaszodi2021}. In fact, a central issue is whether one can use the $2d$ color classes suggested by the volume Helly theorem. The main result of this note is such a Helly type theorem in the plane. The case of higher dimensions remains open.

\begin{thm}\label{th:main} Let $\F_1,\ldots, \F_{4}$ be finite families of convex sets in $\R^2$. If for every transversal $F_1,F_2,F_3,F_4$ of this system $|\bigcap_i^4 F_i|> V$, then $|\bigcap \F_i| > 1$ for some $i\in [4]$. Here $V$ is a universal constant.
\end{thm}

The key step in the proof of this theorem is the following result.

\begin{thm}\label{th:key} Let $Q_1,\ldots, Q_4$ be quadrilaterals in $\R^2$, each of area at most 1. $Q_i$ is the intersection of four fixed halfplanes $H_{i1},H_{i2},H_{i3},H_{i4}$ that form a family $\F_i$ for every $i$. Then there is a transversal of the system $\F_1,\F_2,\F_3,\F_4$ whose intersection has area at most 4096.    
\end{thm}

\section{Some background}\label{sec:back}

We will need the following recent result of~\citet{DeLoera}. 

\begin{thm}\label{th:colorful-doignon}
Let $\mathcal{F}_1, \ldots, \mathcal{F}_{2^d}$ be finite families of convex sets in $\mathbb{R}^d$. If $\bigcap_{i=1}^{2^d} F_i$ contains a lattice point for every transversal $F_i \in \mathcal{F}_i, i\in[2^d]$, then there is an index $k \in[2^d]$ so that $\bigcap \mathcal{F}_k$ contains a lattice point.
\end{thm}

This theorem is a colorful generalization of a classic result of \citet{Doignon1973} which is the special case when all the families are equal, that is, when $\mathcal{F}_1=\ldots=\mathcal{F}_{2^d}$.

The proof of Theorem~\ref{th:key} uses affine unimodular lattices in $\mathbb R^d$. So we define a \defin{full-rank lattice} in $\mathbb R^d$ as the set
\[
        L=\left\{\sum_{i=1}^d m_i v_i : m_i\in\mathbb Z\right\},
\]
where $v_1,\dots,v_d$ are linearly independent vectors in $\mathbb R^d$.

An \defin{affine lattice} $\Lambda$ in $\mathbb R^d$ is a translate of a full-rank
lattice thus it has the form
\[
        \Lambda=x+L,
\]
where $x\in\mathbb R^d$ and $L\subset\mathbb R^d$ is a full-rank lattice.
Equivalently,
\[
        \Lambda=x+A\mathbb Z^d
        =
        \{x+Am:m\in\mathbb Z^d\},
\]
for some $x\in\mathbb R^d$ and some $A\in \mathrm{GL}_d(\mathbb R)$.

Let $\Lambda=x+A\mathbb Z^d$ be an affine lattice.  The \defin{covolume} of $\Lambda$  is the volume of a fundamental parallelepiped of the lattice, so $\operatorname{covol}(\Lambda)= |\det A|.$

{\bf Definition}.
Let $Y_d$ be the space of affine unimodular lattices in $\mathbb R^d$. Thus an element of $Y_d$ is an affine lattice of the form $\Lambda=g\mathbb Z^d+x$, where $g\in\mathrm{SL}(d,\mathbb R),\ x\in\mathbb R^d$. Let $\mu_d$ be the normalized Haar probability measure on $Y_d$.

When we say that $\Lambda$ is a random affine unimodular lattice, we mean that
$\Lambda$ is the identity random variable on the probability space
$(Y_d,\mu_d)$.  
Thus, for any measurable event $E\subseteq Y_d$, we have 
\[
        \mathbb P(\Lambda\in E):=\mu_d(E).
\]
In particular, for a measurable set $A\subset\mathbb R^d$,
we have $\mathbb P(\Lambda\cap A=\varnothing)
        =
        \mu_d\bigl(\{\Gamma\in Y_d:\Gamma\cap A=\varnothing\}\bigr)$.

For the proof of Theorem~\ref{th:key} we are going to use a beautiful and recent result of Athreya~\cite{Athreya2015}.

\begin{thm}\label{thm:Athreya}
Let $A \subset \mathbb{R}^d$ be a measurable set. Then
$$
\mu\left(\Lambda \in Y_d: \Lambda \cap A=\emptyset\right)<\frac{1}{1+|A|}.
$$
\end{thm}

Simple rescaling gives the following variant (we omit the proof).
\begin{lem}\label{lem:random-affine-lattice}
Let $S\subset\mathbb R^2$ be measurable, and let $\Lambda$ be a random affine
lattice of covolume $\tau$. Then
\[
        \mathbb P(\Lambda\cap S=\varnothing)
        <
        \frac{\tau}{\tau+|S|}
\]
whenever $|S|<\infty$. If $|S|=\infty$, then $\mathbb P(\Lambda\cap S=\varnothing)=0$.
\end{lem}

\section{Proof of Theorem~\ref{th:key}}\label{sec:key}

Assume, for contradiction, that every colorful intersection
\[
        K_{\mathbf j}
        :=
        H_{1j_1}\cap H_{2j_2}\cap H_{3j_3}\cap H_{4j_4},
        \qquad
        \mathbf j=(j_1,j_2,j_3,j_4)\in\{1,2,3,4\}^4,
\]
has area strictly larger than $4096$.  Empty and lower-dimensional intersections
have area $0$, so each $K_{\mathbf j}$ is either of finite area $>4096$, or
has infinite area.

Let $U:=Q_1\cup Q_2\cup Q_3\cup Q_4$. Since $|Q_i|\le 1$, we have $|U|\le 4$. We choose a random affine lattice $\Lambda$ of covolume $\tau=8$. We have
\[
        \mathbb P(\Lambda\cap U\ne\varnothing)
        \le
        \mathbb E|\Lambda\cap U|
        =
        \frac{|U|}{8}
        \le
        \frac12 .
\]

For a fixed colorful intersection $K_{\mathbf j}$,~\Cref{lem:random-affine-lattice}
gives
\[
        \mathbb P(\Lambda\cap K_{\mathbf j}=\varnothing)
        <
        \frac{8}{8+4096}
        =
        \frac1{513}
\]
in the finite-area case, and gives probability $0$ in the infinite-area case.
There are $4^4=256$ colorful choices, hence
\[
        \mathbb P\!\left(
        \exists \mathbf j:\ \Lambda\cap K_{\mathbf j}=\varnothing
        \right)
        <
        \frac{256}{513}
        <
        \frac12 .
\]
Therefore, with positive probability, both $\Lambda\cap U=\varnothing$ and $\Lambda\cap K_{\mathbf j}\ne\varnothing$
for every $\mathbf j\in\{1,2,3,4\}^4$ hold.  
Fix such a lattice $\Lambda$.
We next apply~\Cref{th:colorful-doignon} to the four color classes
\[
        \mathcal F_i:=\{H_{i1},H_{i2},H_{i3},H_{i4}\},
        \qquad i=1,2,3,4.
\]
Since every colorful choice intersects $\Lambda$, the theorem gives some
$i\in\{1,2,3,4\}$ such that
\[
        \Lambda\cap\bigcap_{j=1}^4 H_{ij}\ne\varnothing .
\]
But $\bigcap_{j=1}^4 H_{ij}=Q_i\subset U$, contradicting $\Lambda\cap U=\varnothing$.  Hence some colorful intersection has area at most $4096$.

\section{Proof of Theorem~\ref{th:main}}\label{sec:main}

A standard argument (which we delete) shows that it suffices to consider the case when all the sets in every $\F_i$ are halfplanes. We can assume that each $|\F_i|\ge 4$ by repeating a halfplane in $\F_i$ if necessary. 

We claim that for every $i\in[4]$ there is an $\F_i^*\subset \F_i$ with $|\F_i^*|\le 4$ and $|\bigcap \F_i^*|\le v(2)$ where $v(2)$ is the constant in~\Cref{th:volHelly}, case $d=2$. Fix $i\in[4]$ and assume that every fourtuple of halfplanes in $\F_i$ intersect in a set of area larger than $v(2)$. Theorem~\ref{th:volHelly} implies then that $|\bigcap \F_i|>1$, a contradiction. So there is a fourtuple $H_{i1},H_{i2},H_{i3},H_{i4}\in \F_i$ such that the quadrilateral $Q_i=\bigcap_{j=1}^4 H_{ij}$ has area at most $v(2)$. Let $\F_i^* \subset \F_i$ consist of these four halfplanes.

We now have four quadrilaterals $Q_1,\ldots,Q_4$ each of area at most $v(2).$ Theorem~\ref{th:key} shows that the intersection of a transversal of the system $\F_1^*,\F_2^*,\F_3^*,\F_4^*$ has area at most $4096 v(2)=2^{12}v(2)$. \qed

For an explicit value of the constant, we invoke Brazitikos'
quantitative Helly theorem~\cite{Brazitikos2017}. In the planar case his result allows us to take
\[
        v(2)=(2\sqrt[3]{\pi}\cdot 2)^3=64\pi .
\]
Therefore the preceding argument gives the colorful constant
\[
        V=4096\,v(2)=2^{12}\cdot 64\pi=2^{18}\pi .
\]

\bibliographystyle{plainurlnat}
\bibliography{ref}

@incollection{Athreya2015,
  author    = {Athreya, Jayadev S.},
  title     = {Random affine lattices},
  booktitle = {Geometry, groups and dynamics},
  series    = {Contemporary Mathematics},
  volume    = {639},
  pages     = {169--174},
  publisher = {American Mathematical Society},
  year      = {2015},
  doi       = {10.1090/conm/639/12793}
}

@article{Helly1923,
  author  = {Helly, Eduard},
  title   = {{\"U}ber Mengen konvexer K{\"o}rper mit gemeinschaftlichen Punkten},
  journal = {Jahresbericht der Deutschen Mathematiker-Vereinigung},
  volume  = {32},
  year    = {1923},
  pages   = {175--176}
}

@article{BaranyKatchalskiPach1984,
  author  = {B{\'a}r{\'a}ny, Imre and Katchalski, Meir and Pach, J{\'a}nos},
  title   = {Helly's Theorem with Volumes},
  journal = {The American Mathematical Monthly},
  volume  = {91},
  number  = {6},
  year    = {1984},
  pages   = {362--365},
  doi     = {10.1080/00029890.1984.11971432}
}

@article{Barany1982,
  author  = {B{\'a}r{\'a}ny, Imre},
  title   = {A Generalization of {Carath{\'e}odory}'s Theorem},
  journal = {Discrete Mathematics},
  volume  = {40},
  number  = {2--3},
  year    = {1982},
  pages   = {141--152},
  doi     = {10.1016/0012-365X(82)90115-7}
}

@article{DamasdiFoldvariNaszodi2021,
  author  = {Dam{\'a}sdi, G{\'a}bor and F{\"o}ldv{\'a}ri, Vikt{\'o}ria and Nasz{\'o}di, M{\'a}rton},
  title   = {Colorful Helly-type theorems for the volume of intersections of convex bodies},
  journal = {Journal of Combinatorial Theory, Series A},
  volume  = {178},
  year    = {2021},
  pages   = {105361},
  doi     = {10.1016/j.jcta.2020.105361}
}

@article{Brazitikos2017,
  author  = {Brazitikos, Silouanos},
  title   = {Brascamp--Lieb inequality and quantitative versions of Helly's theorem},
  journal = {Mathematika},
  volume  = {63},
  number  = {1},
  pages   = {272--291},
  year    = {2017},
  doi     = {10.1112/S0025579316000255}
}

@article{DeLoera,
  author  = {De Loera, Jes{\'u}s A. and La Haye, Raymond N. and Oliveros, D{\'e}borah and Rold{\'a}n-Pensado, Edgardo},
  title   = {Helly numbers of algebraic subsets of {$\mathbb{R}^d$} and an extension of Doignon's theorem},
  journal = {Advances in Geometry},
  volume  = {17},
  number  = {4},
  year    = {2017},
  pages   = {473--482},
  doi     = {10.1515/advgeom-2017-0028}
}

@article{Doignon1973,
  author  = {Doignon, Jean-Paul},
  title   = {Convexity in cristallographical lattices},
  journal = {Journal of Geometry},
  volume  = {3},
  pages   = {71--85},
  year    = {1973},
  doi     = {10.1007/BF01949705}
}

\end{document}